\documentclass[conference]{IEEEtran}
\IEEEoverridecommandlockouts
\usepackage{microtype}
\usepackage{latex_styles_common}
\usepackage[font=small]{caption}
\usepackage[sort]{cite}

\usepackage[nameinlink]{cleveref}
\Crefname{figure}{Figure}{Figures}
\crefname{figure}{Figure}{Figures}
\crefname{example}{Example}{Example}
\crefname{theorem}{Theorem}{Theorem}
\crefname{corollary}{Corollary}{Corollary}
\crefname{lemma}{Lemma}{Lemma}
\crefname{proposition}{Proposition}{Proposition}
\crefname{assumption}{Assumption}{Assumption}
\crefname{section}{Section}{Section}
\crefname{algorithm}{Algorithm}{Algorithm}

\usepackage{enumitem}
\newlist{propenum}{enumerate}{1} 
\setlist[propenum]{label=\alph*{\rm)}, ref=\theproposition(\alph*)}
\crefalias{propenumi}{proposition}
\newlist{corenum}{enumerate}{1} 
\setlist[corenum]{label=\alph*{\rm)}, ref=\thecorollary(\alph*)}
\crefalias{corenumi}{corollary}
\newlist{lemenum}{enumerate}{1} 
\setlist[lemenum]{label=\alph*{\rm)}, ref=\thelemma(\alph*)}
\crefalias{lemenumi}{lemma}

\usepackage{amsthm,thmtools}
\declaretheorem[name=Theorem,numberwithin=section]{theorem}

\declaretheorem[name=Proposition,numberlike=theorem]{proposition}

\declaretheorem[name=Assumption,numberlike=theorem]{assumption}

\numberwithin{equation}{section}
\numberwithin{theorem}{section}
\numberwithin{figure}{section}
\numberwithin{proposition}{section}
\numberwithin{example}{section}
\numberwithin{definition}{section}
\numberwithin{assumption}{section}
\numberwithin{lemma}{section}

\input{macros}

\newcommand{\convA}{\widehat\Ascr}

\newcommand{\FaceA}{\Fscr_{\!\!\scriptscriptstyle\Ascr}}

\newcommand{\gauge}{\gamma}
\newcommand{\As}{_{\scriptscriptstyle\Ascr}}

\DeclareMathOperator{\supp}{\hbox{\rm\textbf{supp}}}

\begin{document}

\title{Bundle methods for dual atomic pursuit\thanks{Date:
    \today. This work was supported by ONR award N00014-16-1-2242.}  }

\author{
\IEEEauthorblockN{Zhenan Fan}
\IEEEauthorblockA{\textit{Department of Computer Science} \\
\textit{University of British Columbia}\\
Vancouver, Canada\\
zhenanf@cs.ubc.ca}
\and
\IEEEauthorblockN{Yifan Sun}
\IEEEauthorblockA{\textit{Department of Computer Science} \\
\textit{University of British Columbia}\\
Vancouver, Canada\\
ysun13@cs.ubc.ca}
\and
\IEEEauthorblockN{Michael P. Friedlander}
\IEEEauthorblockA{\textit{Department of Computer Science} \\
\textit{University of British Columbia}\\
Vancouver, Canada\\
mpf@cs.ubc.ca}
}

\maketitle

\begin{abstract}
  The aim of structured optimization is to assemble a solution, using
  a given set of (possibly uncountably infinite) atoms, to fit a model
  to data. A two-stage algorithm based on gauge duality and
  bundle method is proposed.  The first stage discovers the
  optimal atomic support for the primal problem by solving a sequence
  of approximations of the dual problem using a bundle-type
  method. The second stage recovers the approximate primal
  solution using the atoms discovered in the first stage. The overall
  approach leads to implementable and efficient algorithms for large
  problems.
\end{abstract}

\section{Introduction}

A recurring approach for solving inverse problems that arise in
statistics, signal processing, and machine learning is based on
recognizing that the desired solution can often be represented as the
superposition of a relatively few canonical atoms as compared to the
signal's ambient dimension. Canonical examples include compressed
sensing and model selection, where the aim is to obtain sparse vector
solutions; and recommender systems, where low-rank matrix solutions
are required. Our aim is to design a set of algorithms that leverages
this sparse atomic structure in order to gain computational
efficiencies necessary for large problems.

Define the set of atoms by a set $\Ascr\subset\Re^n$. The atomic set may
be finite or infinite, but in either case, we assume that the set is
closed and bounded.  A point $x\in\Re^n$ is said to be sparse with
respect to $\Ascr$ if it can be written as a nonnegative superposition
of a few atoms in $\Ascr$, i.e.,
\[
  x = \sum\limits_{a \in \Ascr} c_a a, \quad c_a \geq 0,
\]
where most of the coefficients $c_a$ associated with each atom $a$ are
zero. Two examples are compressed sensing, where the atoms are the
canonical unit vectors and a sparse decomposition is equivalent to sparsity in x
and low-rank matrix completion; where the atoms are the set of
rank-1 matrices and a sparse decomposition is equivalent to low rank. In
each of these cases, there is a convex optimization problem whose
solution is sparse relative to the required atomic set. There now
exists a substantial literature that delineates conditions under which
the correct solution is identified, typically in a probabilistic sense
\cite{recht2008necessary,donoho2006most,candes2004robust,recht2010guaranteed}.

Our focus here is on the approach advocated by~Chandrasekaran et
al.~\cite{chandrasekaran2012convex}, who identified a set of convex
analytical techniques based on gauge functions, which are norm-like
functions that are especially well suited to the atomic description of
the underlying model. We describe below a linear inverse
problem that generalizes the models analyzed by
Chandrasekaran et al.

\section{Atomic pursuit}

The atomic set $\Ascr$ induces the gauge function
\begin{equation}
  \label{eq:gauge}
  \gamma_{\Ascr}(x) = \inf \Set{\mu\ge0 | x \in \mu \hat \Ascr },
\end{equation}
where $\convA = \conv(\Ascr \cup \{0\})$ denotes the convex hull of $\Ascr$ and $0$. 
The gauge to $\Ascr$ can be expressed equivalently as
\begin{equation}
  \label{eq:gauge2}
  \gamma_\Ascr(x) = \inf\bigg\{\sum_{a\in\Ascr, c_a\ge0}c_a \Bigm\vert x=\sum_{a\in\Ascr}c_aa \bigg\};
\end{equation}
see Bonsall~\cite{bonsall1991general}.  The atomic pursuit problem
minimizes the gauge function over a set of linear measurements
$M x\in\Bscr$, where $M:\Re^n\to\Re^m$ is a linear operator, and
$\Bscr\subset\Re^m\setminus\{0\}$ denotes the admissible set:
\begin{equation}
\label{problem:atomic_pursuit}
\minimize{x\in \Re^n} \quad \gamma_{\Ascr}(x) \quad\st\quad M x \in \Bscr.
\end{equation}
Chandrasekaran et al. \cite{chandrasekaran2012convex} and Amelunxen et
al.~\cite{amelunxen2014living} describe conditions under which a
solution to this convex optimization problem yields a good
approximation to the underlying ground truth.

\begin{table*}
\centering
\resizebox{\textwidth}{!}{  
\begin{tabular}{@{}lccccc@{}}    
\toprule 
Atomic sparsity & $\Ascr$ & $\hat \Ascr$ & 
  $\gauge\As(x)$ & $\supp\As(x)$ & $\sigma\As(z)$ \\
\midrule
non-negative & $\cone(\set{e_1,\ldots,e_n})$ & non-negative orthant &
  $\delta_{\geq 0}$ & $\cone(\set{e_i \mid x_i > 0})$ & $\delta_{\leq 0}$\\

elementwise & $\set{ \pm e_1,\ldots,\pm e_n }$ & cross polytope & 
  $\|\cdot\|_1$ & $\set{\sign(x_i)e_i \mid x_i \neq 0}$ & $\|\cdot\|_\infty$\\ 

low rank & $\set{uv^T \mid \|u\|_2=\|v\|_2 = 1}$ & nuclear-norm ball & 
  nuclear norm  & singular vectors of $x$ & spectral norm\\

PSD \& low rank & $\set{uu^T \mid \|u\|_2 = 1}$ & $\set{X\succeq 0 \mid\trace X \leq 1}$ & 
  $\trace + \delta_{\succeq 0}$  & eigenvectors of $x$ & $\max\set{\lambda_{\max}, 0}$\\ 
\bottomrule 
\end{tabular}}
\caption[Commonly used sets of atoms and their gauge and support function
representations]{Commonly used sets of atoms and their gauge and support function
representations. The indicator function $\delta_{\Cscr}(x)$ is zero if $x \in \Cscr$ and $+\infty$  otherwise. 
\label{tab:common-atoms}}  
\end{table*}

Although~\eqref{problem:atomic_pursuit} is convex and in theory
amenable to efficient algorithms, in practice the computational and
memory requirements of general-purpose algorithms are prohibitively
expensive. However, algorithms specially tailored to recognize the
sparse atomic structure can be made to be effective in practice. In
particular, if we had information on which atoms participate
meaningfully in constructing a solution $x^*$,
then~\eqref{problem:atomic_pursuit} can be reduced to a problem over
just those atoms.

Formally, define the set of supports of a vector $x$ with respect to
$\Ascr$ to be all the sets $\Sscr\in \Ascr$ that satisfy
\begin{equation}
  \label{eq:atomic_support}
    x=\sum_{a\in\Sscr}c_a a,\enspace c_a>0, \enspace \gamma_\Ascr(x) = \sum_{a\in \Sscr} c_a,
\end{equation}
i.e., all sets of atoms in $\Ascr$ that contribute non-trivially to the
construction of $x$.  If we can identify any support set
$\Sscr\in\supp_\Ascr(x^*)$ for any solution $x^*$,
then~\eqref{problem:atomic_pursuit} is equivalent to the reduced
problem
\begin{equation}
\label{problem:atomic_pursuit_reduced}
\minimize{x\in \Re^n} \quad \gamma_{\Sscr}(x) \quad \st \quad M x \in \Bscr.
\end{equation}
This is a potentially easier problem to solve, particularly in the
case where it is possible to identify a support set
$\Sscr\in\supp_\Ascr(x^*)$ that has small cardinality. 
For example, when $\convA$ is the cross-polytope, then identifying a
small support $\Sscr$ means that the reduced
problem~\eqref{problem:atomic_pursuit_reduced} only involves the few
variables in $\Sscr$.  Similarly, when $\convA$ is the set of rank-1
positive semidefinite matrices, then identifying a small support $\Sscr$
corresponds to finding the eigenspace of a low-rank solution $X^*$. In
both cases, knowing $\Sscr$ can reduce the computational complexity
significantly.

\section{Dual atomic pursuit}

Our approach for constructing the optimal support set $\supp_\Ascr^*$ is
founded on approximately solving the a dual problem that is particular
to gauge optimization~\eqref{problem:atomic_pursuit}. These dual
problems take the form
\begin{equation}
\label{problem:atomic_pursuit_gauge_dual}
\minimize{y \in \Re^m} \quad \sigma\As(M^* y) \quad \st \quad  y \in \Bscr',
\end{equation}
where $\sigma\As(z) = \sup_{a\in\convA}\,\ip{a}{z}$ is the
support function to the set $\convA$, and
$\Bscr' = \Set{y \in \Re^m| \ip{b}{y} \geq 1\ \forall b \in \Bscr} $ is
the \emph{antipolar} to $\Bscr$. The dual relation between the pair~\eqref{problem:atomic_pursuit} and~\eqref{problem:atomic_pursuit_gauge_dual} is encapsulated by the inequality
\begin{equation}
  \label{eq:strong_duality}
  1\le\ip{x}{M^*y}\le\gamma_\Ascr(x)\cdot\sigma\As(M^*y),
\end{equation}
which holds for all pairs $(x,y)$ that are primal-dual feasible, i.e.,
$Mx\in\Bscr$ and $y\in\Bscr'$. Moreover, under a suitable constraint
qualification, $(x,y)$ is optimal if and only if all of the above
inequalities hold with equality, in which case strong duality holds
\cite[Corollary~5.4]{friedlander2014gauge}.

The following theorem shows that the gauge dual reveals the optimal
support for the primal solution.
\begin{theorem}[Optimal support identification]
  \label{thm:dual-support-id}
  Let $(x, y) \in \Re^n \times \Re^m$ be any optimal primal-dual
  solution of the dual pair~\eqref{problem:atomic_pursuit}
  and~\eqref{problem:atomic_pursuit_gauge_dual}. Then
  \[
    \Sscr\subseteq \Fscr\As(M^*y) \quad \forall \Sscr\in \supp_\Ascr^*,
  \]
  where $\Fscr\As(z) := \set{x\in\Ascr|\ip x z = \sigma\As(z)} = \partial\sigma\As(z)$ is the face of $\convA$ exposed by $z$.
\end{theorem}

This  result can be interpreted geometrically: the optimal support
atoms lie in the face of $\convA$ exposed by $M^*y$. Moreover,
each atom is a subgradient of $\sigma_\Ascr$. The theoretical basis for
this approach is outlined by Friedlander et al.~\cite{friedlander2014gauge} and Aravkin et al.~\cite{aravkin2017foundations}.

\section{Bundle-type two-stage algorithm}

The cutting-plane method for general nonsmooth convex optimization was first
introduced by Kelley~\cite{kelley1960cutting}. It solves the optimization
problem via approximating the objective function by a bundle of linear
inequalities, called cutting planes. The approximation is iteratively refined
by adding new cutting planes computed from the responses of the oracle. The
method is not to approximate the objective function over its entire domain 
by a convex polyhedron, but to construct an approximate valid lower minorant 
near the optimum. Several
stabilized versions, usually known as bundle methods, were subsequently
developed by Lemarechal et al.~\cite{lemarechal1995new} and
Kiewel~\cite{kiwiel1990proximity}.

We give a simplified description of the construction of the lower minorant in the
context of a generic convex function $f: \Re^n \to \Re$. Let $\set{x\j, g\j
\in \partial f(x\j)}_{j = 1}^k$ be the set of pairs of iterates and
subgradients visited through iteration $k$. The cutting plane model at
iteration $k$ is

\begin{equation}
  \label{eq:cutting_plane}
  f\k(x) = \max\limits_{j = 1, \dots, k}\Set{f(x\j) + \ip{g\j}{x - x\j}}.
\end{equation}
In our simplified description, the cutting-plane model is polyhedral.
We can, however, define these more
generally, as described in~\cref{sec:spec}, where the
models are spectrahedral.

The next proposition shows that, when specialized to support functions, the
cutting-plane models are themeselves support functions, which only depends on
the previous subgradients.

\begin{proposition}(Cutting-plane model for support functions)
  \label{thm:inscribing-atoms}
  The cutting-plane model~\eqref{eq:cutting_plane} for $f=\sigma\As$ and $\set{z\j, a\j}_{j = 1}^k$ being the set of pairs of iterates and subgradients is
  \[
    f\k(z) = \sigma_{\Ascr\k}(z)\quad\mbox{with}\quad \Ascr\k = \{a^{(1)},
    \dots, a\k\}.
  \]
\end{proposition}

It follows that the cutting-plane model for this objective of~\eqref{problem:atomic_pursuit_gauge_dual} takes the form
\[
  \sigma_{\Ascr_k}(M^*y) \quad\mbox{with}\quad \Ascr_k\subset\Ascr_{k+1}\subset\cdots\subset\convA.
\]
Coupled with Theorem~\ref{thm:dual-support-id}, we observe that that
the sets $\Ascr_k \subset \convA$ that define the cutting-plane model
are constructed from the faces of $\convA$ exposed by previous
iterates $\{M^*y_i\}_{i = 1}^k$. The sets $\Ascr_k$ thus contain atoms
that are candidates for the support of the optimal solution. In order
to make this approach computationally useful, care must be taken to
ensure that the sets $\Ascr_k$ do not grow too large.
Proposition~\ref{thm:inscribing-atoms} is thus most useful as a guide,
and we consider below an algorithmic variation that allows us to
periodically trim the inscribing sets.

Our method is based on the level bundle method
introduced by Bello Cruz and Oliveira~\cite{cruz2014level}. Each iterate is
computed via a projection onto the level set of the corresponding lower
minorant. The sequence of candidate atomic sets $\set{\Ascr\k}$ inscribes $\convA$,
but does not necessarily grow monotonically as per
Proposition~\ref{thm:inscribing-atoms}. Instead, we follow the recipe outlined
by Br\"annlund and Kiwiel~\cite{brannlund1995descent}, which only requires $\Ascr\kp1$ to contain the elements
that contribute to $y\kp1$. However, this only works for the polyhedral atomic
sets $\Ascr\k$. Applied to more general atomic sets, not necessarily polyhedral, 
\begin{equation}
  \label{eq:convergence_guarantee}
  \Fscr_{\Ascr\k}(z\kp1)\cup \{a\kp1\}
  \subseteq
  \Ascr\kp1
\end{equation}
where $y\kp1$ is the latest iterate, $z\kp1:=M^*y\kp1$, and $a\kp1\in\Fscr\As(z\kp1)$. 
This rule ensures that the updates to the candidate atomic set $\Ascr\k$ always
contain exposed atoms that define the lower minorant, and at least one exposed
atom from the full set. The general version of our proposed method is outlined
in~\cref{alg:general}.

\begin{algorithm}[t]
   \caption{Generic bundle method}
   \label{alg:general}
   \hspace*{\algorithmicindent} \textbf{Input}: \\
   \hspace*{\algorithmicindent} $\delta > 0$ (tolerance)\\ 
   \hspace*{\algorithmicindent} $y^{(1)} \in \Bscr'$ (initial point)\\
   \hspace*{\algorithmicindent} $d^*$ (optimal dual value)
\begin{algorithmic}[1]
   \STATE (Initialize bundle)\quad Construct $\Ascr^{(1)}$ such that \vspace{-2mm}
   \[\Ascr^{(1)} \subseteq \Fscr\As(M^* y^{(1)})\]\vspace{-6mm}
   \STATE (Set center)\quad $\hat y = y^{(1)}$
   \FOR{$k = 1,2,...$}
      \STATE (Upper bound)\quad $U\k = \min\limits_{i = 1, \dots. k} \sigma\As(M^* y^{(i)})$
      \STATE (Gap)\quad $\delta\k = U\k - d^*$
      \STATE (Stopping criterion)\quad \textbf{If} $\delta\k \leq \delta$ \textbf{then stop}
       \STATE (Level set)\vspace{-2mm}
      \begin{align*}
         L\k &:= \Set{y \mid \sigma_{\Ascr\k}(M^*y) \leq d^*}\\
         H\k &:= \Set{y \mid \ip{y - y\k}{\hat y - y\k} \leq 0}\\
         Y\k &:= L\k \cap H\k \cap \Bscr'
       \end{align*} \vspace{-6mm}
      \STATE (Next iterate) \vspace{-2mm}
      \[y\kp1 = \proj_{Y\k}(\hat y)\]\vspace{-6mm}
      \STATE (Update bundle) Construct $\Ascr\kp1$ satisfying~\eqref{eq:convergence_guarantee}
   \ENDFOR
   \RETURN $\Ascr\k$
\end{algorithmic}
\end{algorithm}

Let $\Sscr^*$ denote the solution set to the
problem~\eqref{problem:atomic_pursuit_gauge_dual}, and assume that $\Sscr^*
\neq \emptyset$. Our next theorem shows the convergence of~\cref{alg:general}.

\begin{theorem}[Convergence of~\cref{alg:general}] \label{thm1}

  The sequence $\set{y\k}_{k = 1}^\infty$ converges to the point $y^*:=
  \proj_{\Sscr^*}(\hat y)$. 

\end{theorem} 
The proof for this theorem follows directly from Bello Cruz and Oliveira~\cite[Theorem~3.4]{cruz2014level} with small modification.

We propose two approaches for constructing the candidate sets $\Ascr_k$
specialized for polyhedral atomic sets and for spectral atomic
sets. In the polyhedral case, we construct $\Ascr_k$ by the traditional
polyhedral cutting-plane model described by
Proposition~\ref{thm:inscribing-atoms} and show in that the optimal
atomic support is identified in finite time. In the spectral case, we
follow Helmberg and Rendl~\cite{helmberg2000spectral} who replace the
polyhedral cutting-plane model by a semidefinite cutting-plane
model. Principly, the convergence of both methods follows from
Br\"annlund et al.~\cite[Theorem~3.7]{brannlund1995descent}.

For simplicity, we assume that we know the optimal value $d^*$
of~\eqref{problem:atomic_pursuit_gauge_dual}. This assumption is valid
in cases where it is possible to normalize the ground truth, and then
the optimal value $d^*$ is known in advance.

\subsection{Polyhedral Atomic Set} \label{sec:poly}
In the case where the atomic set $\Ascr$ is finite, and the convex hull
$\convA$ is polyhedral. Our specialized bundle update, which satisfies~\eqref{eq:convergence_guarantee}, is given by
\begin{equation}\label{eq:poly_update}
    \Ascr\kp1 = \Fscr_{\delta, \Ascr}(z\kp1) \cup \{a\kp1\}
\end{equation}
where 
\[\Fscr_{\delta, \Ascr}(z\kp1) = \set{a \in \Ascr\k \mid \ip{a}{z\kp1} \geq \sigma_{\Ascr\k}(z\kp1)- \delta}\]
is the relaxed exposed face. 

Assume the first stage stops at iteration
$T$ with a candidate atomic set $\Ascr^{(T)}$. Then we solve the primal
atomic pursuit on the discovered atoms, namely the reduced
problem~\eqref{problem:atomic_pursuit_reduced} with $\Sscr=\Ascr^{(T)}$.
The following theorem shows the recovery guarantee. 

\begin{theorem}\label{thm:poly_guarantee}
  If the atomic set is finite, then~\cref{alg:general} with specialized bundle update~\eqref{eq:poly_update}  terminates in a finite number of iterations $T$ for all $\delta \geq 0$: 
    \begin{itemize}
      \item if $\delta > 0$, then 
    \[
        0 \le \gauge\As(\overline{x}) - \gauge\As(x^*) \leq \frac{\delta}{d^*(d^*-\delta)},
      \]
      where $x^*$ and $\overline{x}$ denote the optimizer for atomic pursuit and
      recovered solution respectively;
      \item if $\delta = 0$, then $\Sscr \subseteq \Ascr_T$, for all $\Sscr \in \supp\As(x^*)$.
    \end{itemize}
\end{theorem}

\subsection{Spectral Atomic Set} \label{sec:spec}
The bunlde method on the gauge
dual~\eqref{problem:atomic_pursuit_gauge_dual} that we have so far described can be
interpreted as forming inscribing polyhedral approximations to the atomic set.
However, when the atomic set is not polyhedral, which is usually the case when
dealing with semidefinite programs, these polyhedral bundle-types do not
perform that well. Can we form richer non-polyhedral approximations 
to the atomic set? Helmberg and Rendl~\cite{helmberg2000spectral} instead propose a
semidefinite cutting plane model that is formed by restricting the feasible
set to an appropriate face of the semidefinite cone. Here we will apply a
similar idea. 

The spectral atomic set $\Ascr = \set{uu^T | \|u\|_2 = 1}$ contains
uncountably many atoms. The support function corresponding to $\Ascr$ has the
explicit form $\sigma\As(z) = \max\{\lambda_{\max}(z),
0\}$~\cite[Proposition~7.2]{friedlander2014gauge}. 

Now consider problem~\eqref{problem:atomic_pursuit_gauge_dual}. Let $y\k$ be
the current iterate and $P\k$ be an $n$-by-$r$ orthogonal matrix whose range
intersects the leading eigenspace of $M^*y\k$. (In this setting, the adjoint
operator maps $m$-vectors to $n$-by-$n$ symmetric matrices.) Then we can build
a local spectral inner approximation of $\Ascr$ by

\[
  \overline{\Ascr}\k = \set{P\k V {P\k}^T | V \succeq 0,\ \trace(V) \le 1}.
\]
This definition only uses information from the current iterate $y\k$. Now we
consider all the previous iterates $y^{(1)}, \dots, y\km1$.  Following the
aggregation step proposed by Helmberg and Rendl~\cite{helmberg2000spectral},
we aggragate the information from previous iterates into a single matrix $W\k
\in \convA$, and get a richer spectral inner approximation by

\begin{equation}
  \label{eq:spectral_approximation}
  \Ascr\k = \set{\alpha W\k + P\k V{P\k}^T | \alpha + \trace(V) \le 1,\ V \succeq 0}.
\end{equation}
The following result shows that the corresponding lower minorant is
easy to compute.
\begin{proposition}[Spectral cutting-plane model] \label{prop:spectral_cut_plane}
  With $\Ascr\k$ as defined in~\eqref{eq:spectral_approximation}, the
  spectral cutting-plane model is given by
  \begin{equation}
  \label{eq:spectral_support}
  \sigma_{\Ascr\k}(M^*y) = \max\set{0, \lambda_{\max}(T\k), \ip{W\k}{M^*y}}.
\end{equation}
where $T\k = {P\k}^TM^*yP\k$. 
\end{proposition}

The update of the bundle is as follows. Take any matrix $\Wbar :=
\bar\alpha W\k + P\k\Vbar {P\k}^T \in \Ascr\k$ exposed by the latest iterate $M^*y\kp1$,
i.e., $\ip{\Wbar}{M^*y\kp1} = \sigma_{\Ascr\k}(M^*y\kp1)$. Then the important
information of $\Wbar$ is contained in the spectrum spanned by the
eigenvectors of $\Vbar$ associated with maximal eigenvalues. Consider the
eigenvalue decomposition $\Vbar = Q\Lambda Q^T$, where $\Lambda =
\diag(\lambda_1, \dots, \lambda_r)$ with $\lambda_1 \geq \dots \geq
\lambda_r$. Split the spectrum $\Lambda$ into two parts: $\Lambda_1=\lambda_1
I$ contains the maximal eigenvalue with possible multiplicity, and $\Lambda_2$
contains the remaining eigenvalues. Let $Q_1$ and $Q_2$ be the corresponding
eigenvectors.  Then we update the bundle by

\begin{equation} \label{eq:spectral_bundle_update}
\begin{split}
W\kp1 &= \dfrac{\bar \alpha W\k + P\k Q_2\Lambda_2Q_2^T{P\k}^T}{\bar \alpha + \trace(\Lambda_2)},\\
P\kp1 &= \orthog[P\k Q_1, v\kp1],
\end{split}
\end{equation}
where $v\kp1$ is any leading normalized eigenvector of
$M^*y\kp1$. 

Because the spectral atomic set is a continuum, we do not expect to
exactly obtain the optimal atomic support, and thus exact recovery of
the primal solution in finite time is not possible. Friedlander and
Mac\^edo\cite[Corollary~4]{friedlander2016low} show that an
approximate primal solution can be recovered by solving a semidefinite
least-squares problem. Given a set of candidate optimal atoms $\Ascr^{(T)}$,
an approximate primal solution can be obtained as the solution of 
\begin{equation}
\label{problem:primal_recover_spectral}
\minimize{x\in \Re^n} \enspace \tfrac{1}{2}\|M x - b\|^2_2 \enspace \st \enspace x \in \cone\big(\Ascr^{(T)}\big).
\end{equation}

It is shown by Ding in~\cite{ding2019optimal} that under certain assumptions, the recovery quality can be guaranteed. 
\begin{assumption}[Assumptions to ensure recovery quality] The following three assumptions are critical for the recovery guarantee of~\eqref{problem:primal_recover_spectral}.
\label{assumption} 
  \begin{propenum}
  \item (Uniqueness) Both primal problem~\eqref{problem:atomic_pursuit} and dual problem~\eqref{problem:atomic_pursuit_gauge_dual} have unique solution $x^*$ and $y^*$.
  \item (Strong duality) Every solution pair $(x^*, y^*)$ satisfies strong duality~\cite[Corollary~5.4]{friedlander2014gauge}:
  \[1 = \ip{x^*}{M^*y^*} = \gauge\As(x^*)\sigma\As(M^*y^*).\]
  \item (Strict complementarity) Every solution pair $(x^*, y^*)$ satisfies the strict complementarity condition:
  \[\rank(x^*) = \mbox{ Mutiplicity of}\enspace\lambda_{\max}(M^*y).\]
  \end{propenum}
\end{assumption}

\begin{theorem}[Spectral recovery guarantee]
  Assume~\cref{assumption} holds and the first stage stops at iteration $T$ with
  $U_T - d^* \leq \delta$. Let $x^*$ and $\overline{x}$ respectively denote the
  optimizer for \eqref{problem:atomic_pursuit} and
  \eqref{problem:primal_recover_spectral}. Then
  \[\|x^* - \overline{x}\|_F = \mathcal{O}(\sqrt{\delta})\] for any
  solution $x^*$ of~\eqref{problem:atomic_pursuit}.
\end{theorem}

The proof for this theorem follows directly from Ding et al.~
\cite[Theorem~1.2]{ding2019optimal} with small modification.

\section{Experiments}\label{sec:experiments}
\subsection{Basis pursuit denoising} \label{subsec:bpdn}
The basis pursuit denoising (BPDN) \cite{chen2001atomic} problem
arises in sparse recovery applications. Let $M: \Re^n \to \Re^m$ be some
measurement matrix. Let $x_0$ denote the original signal and
$b = Mx_0$ be the vector ofobservations, where $x_0$ is sparse and the
observation $b$ might be noisy. For some expected noise level
$\epsilon>0$, the BPDN model is
\begin{equation}
\label{problem:basis_pursuit}
\minimize{x\in \Re^n} \enspace \|x\|_1 \enspace \st \enspace \|Mx - b\|_2 \leq \epsilon.
\end{equation}
In this case, the atomic set $\Ascr$ is the set of signed unit one-hot
vectors $\Ascr = \set{\pm e_1,\ldots,\pm e_n}$ and $\Bscr$ is the
2-norm ball centered at $b$ with radius $\epsilon$. The corresponding gauge
dual problem is given by
\begin{equation}
\label{problem:basis_pursuit_dual}
\minimize{y\in \Re^m} \enspace \|M^Ty\|_{\infty} \enspace \st \enspace y \in \Bscr',
\end{equation}
where the antipolar $\Bscr' = \set{y\mid\ip{b}{y} - \epsilon\|y\|_2 \geq 1}$ follows directly from the definition.

\begin{figure}[H]
  \centering
  \includegraphics[width = \linewidth]{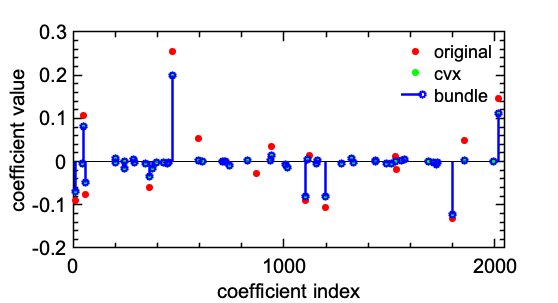}   
  \caption{The coefficients of the original signal and the signals recovered by cvx and by bundle. \label{fig:bp_figure}}
\end{figure}
\begin{figure}[H]
  \centering
  \includegraphics[width = \linewidth]{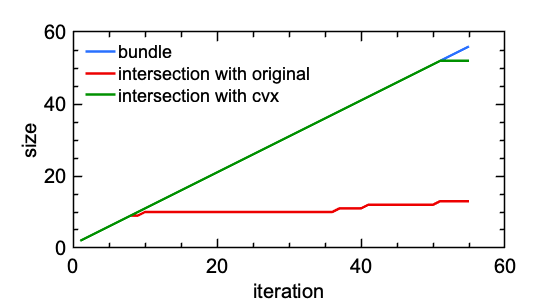}
  \caption{Size of bundle and its intersection with the atomic support for both original signal and signal recovered by cvx. \label{fig:bp_size}
}
\end{figure}

\subsection{Phase retrieval} \label{subsec:phase}
Phase retrieval is a problem of recovering signal from magnitude-only
measurements. Specifically, let $x_0 \in \Re^n$ be some unknown signal
and the measurements are given by $b_i = |\langle x_0, m_i\rangle|^2$,
where each vector $m_i$ encodes illumination $i=1,\ldots,m$. Cand\'es
et al. \cite{candes2015phase} advocate ``lifting'' the signal as
$X_0 = x_0 x_0^T$ so that the measurements are linear in $X_0$:
\[b_i = \ip{x_0x_0^T}{m_im_i^T} = \ip{X_0}{M_i},
  \; i = 1, \dots, k,\] where $M_i = m_im_i^T$. The following
semidefinite program can be used to recover $X^* \approx x_0x_0^T$:
\begin{equation}
\label{problem:phase_retrieval}
\begin{array}{ll}
\minimize{X} & \trace(X) + \delta_{\succeq 0}(X) \\ \st &  \ip{X}{M_i} = b_i, \; i = 1, \dots, m.
\end{array}
\end{equation}
Define $\Ascr$ as the set of normalized positive semidefinite rank-1 matrices, and define the linear operator $M$ as
\begin{equation}
M(X) = \left[ \ip{X}{M_i} \right]_{i=1,\ldots,m}.
\end{equation}
Then the atomic pursuit problem \eqref{problem:atomic_pursuit} is equivalent to \eqref{problem:phase_retrieval} with $\Bscr = \{b\}$. The corresponding gauge
dual problem is
\begin{equation}
\label{problem:phase_retrieval_dual}
\minimize{y\in \Re^m} \enspace \max\set{0, \lambda_{\max}(M^*y)} \enspace \st \enspace \ip{b}{y} \geq 1,
\end{equation}
where the adjoint operator $M^*$ applied to a vector $y$ is defined as
\[M^*y = \sum\limits_{i = 1}^m y_iM_i.\]


\begin{figure}[H]
\centering
  \includegraphics[width=\linewidth]{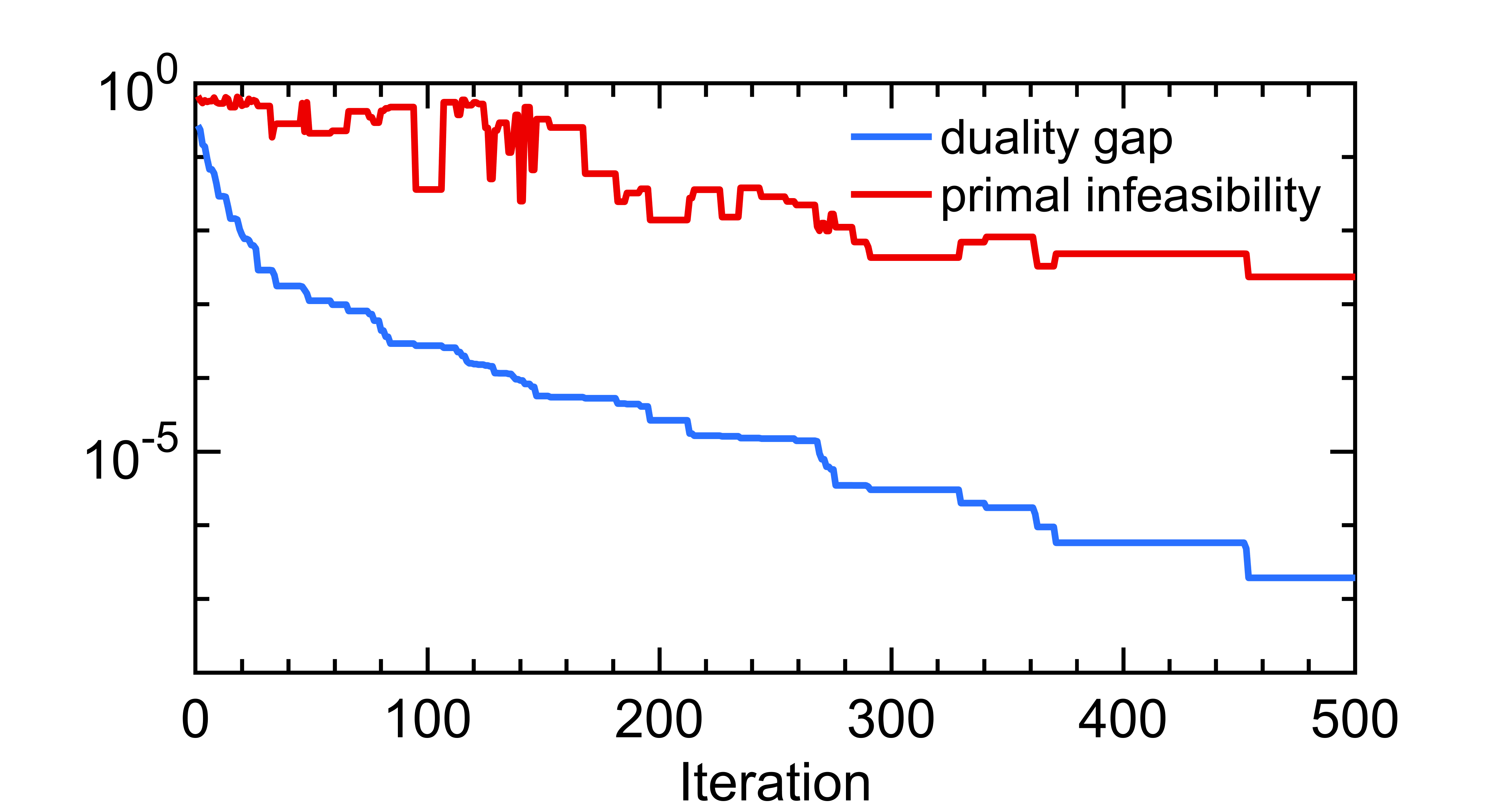}
  \caption{Convergence of objective value and primal infeasibility.}
\label{fig:pl_con}
\end{figure}

\begin{figure}[H]
\centering
  \includegraphics[width=\linewidth]{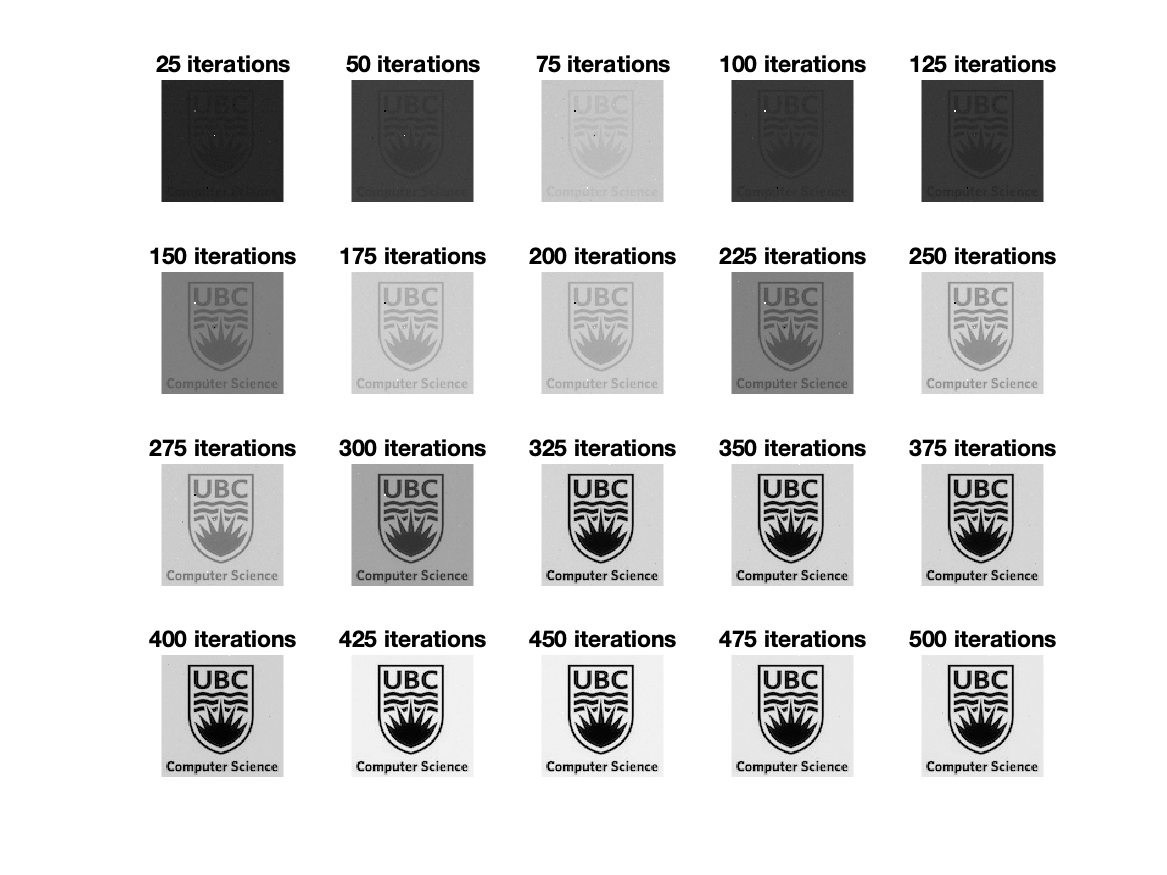}
  \caption{Images recovered at intermediate iterations of the spectral
    bundle method for recovering a
    ground-truth image.
    \label{fig:pl}}
\end{figure}

\section{Conclusion} \label{sec:conclusion}

Convex optimization formulations of inverse problems often come with
very strong recovery guarantees, but the formulations may be too large
to be practical for large problems. This is especially true of
spectral problems, which require very expensive computational
kernels. Our atomic pursuit approach shifts the focus from the
solution of the full convex optimization problem to a sequence of
``reduced'' problems meant to expose the constituent atoms that form
the final solution. In some sense, atomic pursuit is a generalization
of more classical active-set algorithms. Future avenues of research
include the design of specialized SDP solvers for the solution of the
highly-structured bundle subproblems, and applying the algorithm
framework to other atomic sets.

\newpage

\bibliographystyle{IEEEtran} \bibliography{IEEEabrv,Refs}

\begin{thebibliography}{10}
\providecommand{\url}[1]{#1}
\csname url@samestyle\endcsname
\providecommand{\newblock}{\relax}
\providecommand{\bibinfo}[2]{#2}
\providecommand{\BIBentrySTDinterwordspacing}{\spaceskip=0pt\relax}
\providecommand{\BIBentryALTinterwordstretchfactor}{4}
\providecommand{\BIBentryALTinterwordspacing}{\spaceskip=\fontdimen2\font plus
\BIBentryALTinterwordstretchfactor\fontdimen3\font minus
  \fontdimen4\font\relax}
\providecommand{\BIBforeignlanguage}[2]{{%
\expandafter\ifx\csname l@#1\endcsname\relax
\typeout{** WARNING: IEEEtran.bst: No hyphenation pattern has been}%
\typeout{** loaded for the language `#1'. Using the pattern for}%
\typeout{** the default language instead.}%
\else
\language=\csname l@#1\endcsname
\fi
#2}}
\providecommand{\BIBdecl}{\relax}
\BIBdecl

\bibitem{recht2008necessary}
B.~Recht, W.~Xu, and B.~Hassibi, ``Necessary and sufficient conditions for
  success of the nuclear norm heuristic for rank minimization,'' in \emph{2008
  47th IEEE Conference on Decision and Control}.\hskip 1em plus 0.5em minus
  0.4em\relax IEEE, 2008, pp. 3065--3070.

\bibitem{donoho2006most}
D.~Donoho, ``For most large underdetermined systems of linear equations the
  minimal $\ell^1$-norm solution is also the sparsest solution,''
  \emph{Communications on pure and applied mathematics}, vol.~59, no.~6, pp.
  797--829, 2006.

\bibitem{candes2004robust}
E.~Candes, J.~Romberg, and T.~Tao, ``Robust uncertainty principles: Exact
  signal reconstruction from highly incomplete frequency information,''
  \emph{arXiv preprint math/0409186}, 2004.

\bibitem{recht2010guaranteed}
B.~Recht, M.~Fazel, and P.~A. Parrilo, ``Guaranteed minimum-rank solutions of
  linear matrix equations via nuclear norm minimization,'' \emph{SIAM review},
  vol.~52, no.~3, pp. 471--501, 2010.

\bibitem{chandrasekaran2012convex}
V.~Chandrasekaran, B.~Recht, P.~A. Parrilo, and A.~S. Willsky, ``The convex
  geometry of linear inverse problems,'' \emph{Foundations of Computational
  mathematics}, vol.~12, no.~6, pp. 805--849, 2012.

\bibitem{bonsall1991general}
F.~F. Bonsall, ``A general atomic decomposition theorem and banach's closed
  range theorem,'' \emph{The Quarterly Journal of Mathematics}, vol.~42, no.~1,
  pp. 9--14, 1991.

\bibitem{amelunxen2014living}
D.~Amelunxen, M.~Lotz, M.~B. McCoy, and J.~A. Tropp, ``Living on the edge:
  Phase transitions in convex programs with random data,'' \emph{Information
  and Inference: A Journal of the IMA}, vol.~3, no.~3, pp. 224--294, 2014.

\bibitem{friedlander2014gauge}
M.~P. Friedlander, I.~Macedo, and T.~K. Pong, ``Gauge optimization and
  duality,'' \emph{SIAM Journal on Optimization}, vol.~24, no.~4, pp.
  1999--2022, 2014.

\bibitem{aravkin2017foundations}
A.~Y. Aravkin, J.~V. Burke, D.~Drusvyatskiy, M.~P. Friedlander, and K.~MacPhee,
  ``Foundations of gauge and perspective duality,'' \emph{arXiv preprint
  arXiv:1702.08649}, 2017.

\bibitem{kelley1960cutting}
J.~E. Kelley, Jr, ``The cutting-plane method for solving convex programs,''
  \emph{Journal of the society for Industrial and Applied Mathematics}, vol.~8,
  no.~4, pp. 703--712, 1960.

\bibitem{lemarechal1995new}
C.~Lemar{\'e}chal, A.~Nemirovskii, and Y.~Nesterov, ``New variants of bundle
  methods,'' \emph{Mathematical programming}, vol.~69, no. 1-3, pp. 111--147,
  1995.

\bibitem{kiwiel1990proximity}
K.~C. Kiwiel, ``Proximity control in bundle methods for convex
  nondifferentiable minimization,'' \emph{Mathematical programming}, vol.~46,
  no. 1-3, pp. 105--122, 1990.

\bibitem{cruz2014level}
J.~B. Cruz and W.~de~Oliveira, ``Level bundle-like algorithms for convex
  optimization,'' \emph{Journal of Global Optimization}, vol.~59, no.~4, pp.
  787--809, 2014.

\bibitem{brannlund1995descent}
U.~Br{\"a}nnlund, K.~C. Kiwiel, and P.~O. Lindberg, ``A descent proximal level
  bundle method for convex nondifferentiable optimization,'' \emph{Operations
  Research Letters}, vol.~17, no.~3, pp. 121--126, 1995.

\bibitem{helmberg2000spectral}
C.~Helmberg and F.~Rendl, ``A spectral bundle method for semidefinite
  programming,'' \emph{SIAM Journal on Optimization}, vol.~10, no.~3, pp.
  673--696, 2000.

\bibitem{friedlander2016low}
M.~P. Friedlander and I.~Macedo, ``Low-rank spectral optimization via gauge
  duality,'' \emph{SIAM Journal on Scientific Computing}, vol.~38, no.~3, pp.
  A1616--A1638, 2016.

\bibitem{ding2019optimal}
L.~Ding, A.~Yurtsever, V.~Cevher, J.~A. Tropp, and M.~Udell, ``An
  optimal-storage approach to semidefinite programming using approximate
  complementarity,'' \emph{arXiv preprint arXiv:1902.03373}, 2019.

\bibitem{chen2001atomic}
S.~S. Chen, D.~L. Donoho, and M.~A. Saunders, ``Atomic decomposition by basis
  pursuit,'' \emph{SIAM review}, vol.~43, no.~1, pp. 129--159, 2001.

\bibitem{candes2015phase}
E.~J. Candes, Y.~C. Eldar, T.~Strohmer, and V.~Voroninski, ``Phase retrieval
  via matrix completion,'' \emph{SIAM review}, vol.~57, no.~2, pp. 225--251,
  2015.

\end{thebibliography}
\newpage
\begin{appendices}
\section{Proofs}
\subsection{Proof for~\cref{thm:dual-support-id}}
\begin{IEEEproof}
  Let $\Sscr\in\supp\As(x)$. It follows from \eqref{eq:gauge2} that there exist strictly positive
  numbers $\set{c_a|a\in\Sscr}$ such that
  \[
    x^* = \sum_{a\in\Sscr}c_a a \quad\mbox{and}\quad \gauge\As(x) = \sum_{a\in\Sscr}c_a,
  \]
  Let $\hat x = x / \gauge\As(x)$ be a normalized solution. Then
  \[
    \hat x = \sum_{a\in\Sscr}\frac{c_a}{\gauge\As(x)}a
    \quad\mbox{and}\quad
    \gauge\As(\hat x) = \sum_{a\in\Ascr}\frac{c_a}{\gauge\As(x)} \equiv 1.
  \]
  This implies that $\hat x$ is necessarily a strict convex combination of
  every point in $\Sscr$. Thus in order to establish that
  $\Sscr\subseteq\FaceA(M^*y)$, it is sufficient to show
  $\hat x\in\FaceA(M^*y)$. By strong duality,
  \[\ip{\hat x}{M^*y} = \sigma\As(M^*y),\]
  and it follows from the definition of exposed faces that $\hat x \in \FaceA(M^*y^*)$.  
\end{IEEEproof}

\subsection{Proof for~\cref{thm:inscribing-atoms}}
\begin{IEEEproof}
  By the definition of subdifferential of support functions, we have
  \begin{eqnarray*}
    f\k(z) &=& \max\limits_{j = 1, \dots, k} \sigma\As(z\j) + \ip{a\j}{z - z\j}\\
    &=& \max\limits_{j = 1, \dots, k} \ip{a\j}{z}\\
    &=& \sigma_{\Ascr\k}(z).
  \end{eqnarray*}
\end{IEEEproof}

\subsection{Proof for~\cref{thm:inscribing-atoms}}
\begin{IEEEproof}
  By the definition of subdifferential of support functions, we have
  \begin{eqnarray*}
    f\k(z) &=& \max\limits_{j = 1, \dots, k} \sigma\As(z\j) + \ip{a\j}{z - z\j}\\
    &=& \max\limits_{j = 1, \dots, k} \ip{a\j}{z}\\
    &=& \sigma_{\Ascr\k}(z).
  \end{eqnarray*}
\end{IEEEproof}

\subsection{Proof for~\cref{thm:poly_guarantee}}
\begin{IEEEproof}
  \begin{itemize}
    \item When $\delta > 0$, by strong duality(~\cref{prop-support-properties-strong-duality}), we know that 
  \begin{align*}
    \sigma\As(M^* y^*)\gauge\As(x^*) &= 1,\\
    \sigma_{\hat{\Ascr}}(M^* \hat{y})\gauge_{\hat{\Ascr}}(\hat{x}) &= 1.
  \end{align*}
  Then it follows that, 
  \begin{eqnarray*}
    \gauge_{\hat{\Ascr}}(\hat{x}) - \gauge\As(x^*) &=& \frac{1}{\sigma_{\hat{\Ascr}}(M^* \hat{y})} - \frac{1}{\sigma\As(M^* y^*)}\\
    &=& \frac{\sigma\As(M^* y^*) - \sigma_{\hat{\Ascr}}(M^* \hat{y})}{\sigma_{\hat{\Ascr}}(M^* \hat{y})\sigma\As(M^* y^*)}\\
    &\leq& \frac{\epsilon}{d^*(d^* - \epsilon)}.
  \end{eqnarray*}
  Now by the fact that $\gauge_{\hat{\Ascr}}(\hat{x}) \geq \gauge\As(\hat{x}) \geq \gauge\As(x^*)$, we can conclude that 
  \[\gauge\As(\hat{x}) - \gauge\As(x^*) \leq \dfrac{\epsilon}{d^*(d^* - \epsilon)}.\]
  \item When $\delta = 0$, first, we show that the algorithm will terminate in finite steps with stopping criteria~2. Define
  \[S\k_\delta = \{a \in \Ascr \mid \ip{a}{M^*y_{k}} \ge d^* - \delta\}.\]
  By~\cref{thm1}, we know that $y\k \to y^*$, where $y^*$ is some optimal solution to~\eqref{problem:atomic_pursuit_gauge_dual}. Then there exist positive number $K$ such that $\forall k \ge K$, 
  \[\FaceA(M^*y^*) \subseteq S\k_\delta.\]
  And by the construction of $\Ascr\kp1$, we know that $\Ascr\kp1 \subseteq S\k_\delta$ for all $k$. We can thus conclude that there exist some finite number $T$ such that $\Ascr_{T+1} = \Ascr_T$.

  Next, we show that when the algorithm terminate, the bundle will contain $\supp\As(x^*)$. From the discussion above, we know that
  \[\Ascr_{T} = S^{(T)}_\delta\quad\mbox{and}\quad\FaceA(M^*y^*) \subseteq S^{(T)}_\delta.\]
  Then by~\cref{thm:dual-support-id}, the result follows. 
  \end{itemize}
\end{IEEEproof}

\subsection{Proof for~\cref{prop:spectral_cut_plane}}
\begin{IEEEproof}
  \begin{eqnarray*}
    & & \sigma_{\Ascr\k}(M^*y) \\
    &=& \max\bigg\{\ip{M^*y}{\alpha W\k + P\k V{P\k}^T} \mid \\
    & & \quad \alpha \geq 0, \alpha + \trace(V) \le 1, V \succeq 0\bigg\}\\
    &=& \max\limits_{0\le\alpha\le1}\bigg\{\alpha\ip{W\k}{M^*y} + (1 - \alpha)\\
    & & \quad \max\{\ip{P\k V{P\k}^T}{M^*y}\mid \trace(V) \le 1, V \succeq 0\}\bigg\}\\
    &=& \max\limits_{0\le\alpha\le1}\bigg\{\alpha\ip{W\k}{M^*y} + \\
    & & \quad (1 - \alpha)\max\set{0, \lambda_{\max}({P\k}^TM^*yP\k)}\bigg\}\\
    &=& \max\set{0, \lambda_{\max}({P\k}^TM^*yP_k), \ip{W\k}{M^*y}}
  \end{eqnarray*}
\end{IEEEproof}

\end{appendices}
\end{document}